\documentclass[12pt]{article} 
\def\Bbb{\bf} 

\newcommand\C{{ \Bbb C}}

\newcommand\PP{{\tt P}}

\newtheorem{thm}{Theorem}

\newtheorem{cor}{Corollary}
\newtheorem{rmk}{Remark}

\newtheorem{prp}{Proposition}

\def\comment#1{ }
\begin{document}
\title{
The uniformizing differential equation of the complex hyperbolic
structure on the moduli space of marked cubic surface II}
\author{Takeshi SASAKI \thanks{Department of Mathematics, Kobe
University, Kobe 657 Japan} \and Masaaki YOSHIDA \thanks{Department of
Mathematics, Kyushu University 01, Fukuoka 810 Japan}}
\date{May 26, 2000} 
\maketitle
\noindent
{\textbf{Abstract:}}\quad
The relation between the uniformizing equation of the 
complex hyperbolic structure on the moduli space of marked cubic surfaces 
and an Appell-Lauricella hypergeometric system in nine variables is clarified.
\section{Introduction}
In the previous paper (\cite{SaYo2}) we found the uniformizing equation 
which governs the developing map of the 
complex hyperbolic structure on the (4-dimensional) moduli space of marked 
cubic surfaces. Our equation is invariant under the action of the Weyl 
group of type $E_6$. Terasoma and Matsumoto are establishing a theory 
implying that the variation of Hodge 
structure of the cubic surfaces is essentially equivalent to that of a 
certain family of curves which are cyclic covers of the line branching at 
twelve points (\cite{MT}). This theory suggests that  
our uniformizing equation should be equivalent to a 
restriction  of the Appell-Lauricella hypergeometric system $F_D$ in nine variables.
In this paper, we carry out the computation and prove this observation.

\section{Configuration spaces $X(2,N)$}
Let $X(2,N)$ be the configuration space of $N$ (colored) points on the projective line $\PP^1$ defined as
$$X(2,N)=PGL_2\backslash\{(x_1,\dots,x_{N})\in (\PP^1)^N\mid
x_i\not=x_j\ (i\not=j)\}.$$
By normalizing three points as $0, 1, \infty$, the space $X(2,N)$ can be identified with the open affine set
$$\prod_{j=1}^nx_j(1-x_j)\prod_{1\le i<j\le n}(x_i-x_j)\not=0$$
in the affine space coordinatized by $(x_1,\dots,x_n)$, where $n+3=N$.
On this configuration space lives the Appell-Lauricella hypergeometric
system $E_D$, which we review in the next section.
When $n=1$, $X(2,4)$ is isomorphic to $\C-\{0,1\}$.

\section{Appell-Lauricella hypergeometric system $E_D$}
For coordinates $x_1,\dots,x_n$, put
$D_i:=x_i{\partial/\partial x_i}.$
The system
$$D_iD_ju=\sum_{k=1}^np^k_{ij}D_ku + p^0_{ij}u \quad (i,j=1,\dots,n),$$
with parameters $(a,b,c)=(a, b_1,\dots,b_n,c)$, where
\begin{eqnarray*}
p^k_{ij}&=&p^k_{ji},\quad p^0_{ij}=p^0_{ji},\quad 1\le i,j,k\le n;\quad
p^0_{ij}p^k_{ij} =0,\quad i\not=j\not=k\not=i\\
p^i_{ij}&=&b_j {x_j\over x_i-x_j},\quad i\not=j,\qquad
p^0_{ii}=a b_i {x_i\over 1-x_i},\\
p^k_{ii}&=&b_i \left({x_i\over 1-x_i} - {x_i\over x_k-x_i}\right),\quad i\not=k\\
p^i_{ii}&=&-\sum_{k\not=i}b_k {x_k\over x_i-x_k} + 
{(a+b_i)x_i - c+1\over 1-x_i}
\end{eqnarray*}
is called the Appell-Lauricella hypergeometric system of type $D$ and is
denoted by $E^n_D(a, b_1,\dots,b_n,c)$. The Appell-Lauricella
hypergeometric series 
\begin{eqnarray*}
\lefteqn{F_D(a;b_1,\dots,b_n;c|x^1,\dots,x^n)=}\\
&&\sum_{m_1,\dots,m_n=0}^{\infty}{(a,m_1+\cdots+m_n)(b_1,m_1)
\cdots(b_n,m_n)
\over(c,m_1+\cdots+m_n)m_1!\cdots m_n!}(x^1)^{m_1}\cdots(x^n)^{m_n}, \\
\end{eqnarray*}
where $(a,m)=a(a+1)\cdots(a+m-1)$, solves this system around
the origin (\cite{Yo2},\cite{HSY}).
This system has singularities along the divisor defined by
$$\prod_{j=1}^nx_j(1-x_j)\prod_{1\le i<j\le n}(x_i-x_j)=0$$
and at infinity, that is, this system is regular in the configuration space
$X(2,n+3)$. The rank (the dimension of the space of local solutions
 at a(ny) regular point) is $n+1$. When $n=1$, $E_D^1(a,b,c)$ is the Gauss 
hypergeometric equation. 
\par\medskip
Note in general that if a system of the 
above form is of rank $n+1$ then the coefficients $p_{ij}^0$ can be expressed 
in terms of $p_{pq}^r\ (1\le p,q,r\le n)$ and their derivatives; so we often 
describe a system by presenting $p_{ij}^k$ only.

\section{The pull-back of $E_D$ under an embedding of $X(2,n+3)$ into $X(2,2n+4)$}
Let us embed (a Zariski open subset of) the space $X(2,n+3)$ into $X(2,2n+4)$ as
\begin{eqnarray*}
\iota:X(2,n+3)&\ni& (0,\infty,1,x_1,\dots,x_n)\\
&\longmapsto& (0,\infty, 1, x_1,\dots,x_n,-x_1,\dots,-x_n, -1)\in X(2,2n+4).
\end{eqnarray*}
The space $X(2,2n+4)$ is of dimension $2n+1$ and the Appell-Lauricella
hypergeometric system $E_D^{2n+1}$ is of rank $2n+2$.

For a solution $u=u(x_1,\dots,x_{2n+1})$ of $E_D^{2n+1}=E_D^{2n+1}(a,b_1,\dots,b_{2n+1},c)$, put
$$v(x_1,\dots,x_n):=u(x_1,\dots,x_n,-x_1,\dots,-x_n,-1).$$
For generic parameters $(a,b,c)$, the system satisfied by the pull-back $v$ of 
$u$ under $\iota$ would be of rank $2n+2$. We study whether the rank of this system may be of rank $n+1$. 
We have
\begin{eqnarray*}
D_iv&=&(D_i+D_{i+n})u,\quad 1\le i\le n\\
D_iD_jv&=&(D_iD_j + D_iD_{j+n} + D_{i+n}D_j + D_{i+n}D_{j+n})u\\
&=&\sum_{k=0}^{2n+1}(p^k_{i,j}+ p^k_{i,j+n} + p^k_{i+n,j} + p^k_{i+n,j+n})D_ku,
\end{eqnarray*}
where $D_0u:=u$, and the right hand-sides denote the values restricted to 
$$x_{n+1}=-x_1,\dots,\ x_{2n}=-x_n,\quad x_{2n+1}=-1.$$
For $1\le i\not=j\le n$, we have 
\begin{eqnarray*}
D_iD_jv
&=&\sum_{k=1}^{2n}(p^k_{i,j}+ p^k_{i,j+n} + p^k_{i+n,j} + p^k_{i+n,j+n})D_ku,\\
&=&(p^i_{i,j}+p^i_{i,j+n})D_iu + (p^{i+n}_{i+n,j}+p^{i+n}_{i+n,j+n})D_{i+n}u\\
&+&(p^j_{i,j}+p^j_{i+n,j})D_ju + (p^{j+n}_{i,j+n}+p^{j+n}_{i+n,j+n})D_{j+n}u,
\end{eqnarray*}
and
\begin{eqnarray*}
p^i_{i,j}+p^i_{i,j+n}
&=&b_j{x_j\over x_i-x_j} + b_{j+n}{x_{j+n}\over x_i-x_{j+n}}\\
&=&b_j{x_j\over x_i-x_j} + b_{j+n}{-x_{j}\over x_i+x_{j}},\\
p^{i+n}_{i+n,j}+p^{i+n}_{i+n,j+n}
&=&b_j{x_j\over x_{i+n}-x_j} + b_{j+n}{x_{j+n}\over x_{i+n}-x_{j+n}}\\
&=&b_j{x_j\over -x_i-x_j} + b_{j+n}{-x_{j}\over -x_i+x_{j}}.
\end{eqnarray*}
When $i=j$, we have
\begin{eqnarray*}
D_i^2v&=&\sum_{k=0}^{2n+1}(p^k_{i,i}+ 2p^k_{i,i+n} + p^k_{i+n,i+n})D_ku\\
&=&(p^i_{i,i}+2p^i_{i,i+n}+p^i_{i+n,i+n})D_iu
+(p^{i+n}_{i,i}+2p^{i+n}_{i,i+n}+p^{i+n}_{i+n,i+n})D_{i+n}u\\
&+&\sum_{1\le k\not=i\le n}\left[(p^k_{i,i}+p^k_{i+n,i+n})D_ku
+(p^{k+n}_{i,i}+p^{k+n}_{i+n,i+n})D_{k+n}u\right]\\
&+&(p^{2n+1}_{i,i}+p^{2n+1}_{i+n,i+n})D_{2n+1}u + (p^{0}_{i,i}+p^{0}_{i+n,i+n})u,
\end{eqnarray*}
and
\begin{eqnarray*}
p^{2n+1}_{i,i}&=&b_i\left({x_i\over 1-x_i}-{x_i\over x_{2n+1}-x_i}\right)
=b_i\left({x_i\over 1-x_i}-{x_i\over -1-x_i}\right),\\
p^{2n+1}_{i+n,i+n}&=&b_{i+n}\left({x_{i+n}\over 1-x_{i+n}}-{x_{i+n}\over x_{2n+1}-x_{i+n}}\right)
=b_{i+n}\left({-x_i\over 1+x_i}-{-x_i\over -1+x_i}\right),
\end{eqnarray*}
\begin{eqnarray*}
p^k_{i,i}+p^k_{i+n,i+n}
&=&b_i\left({x_i\over 1-x_i}-{x_i\over x_k-x_i}\right)+
b_{i+n}\left({-x_i\over 1+x_i}-{-x_i\over x_k+x_i}\right),\\
p^{k+n}_{i,i}+p^{k+n}_{i+n,i+n}
&=&b_i\left({x_i\over 1-x_i}-{x_i\over -x_k-x_i}\right)+
b_{i+n}\left({-x_i\over 1+x_i}-{-x_i\over -x_k+x_i}\right),
\end{eqnarray*}
\begin{eqnarray*}
p^i_{i,i}&+&2p^i_{i,i+n}+p^i_{i+n,i+n}=
-\sum_{1\le k\not=i\le 2n+1}b_k{x_k\over x_i-x_k}
+{1\over 1-x_i}\left\{(a+b_i)x_i-c+1\right\}\\
&+& 2b_{i+n}{x_{i+n}\over x_i-x_{i+n}}
+b_{i+n}\left({x_{i+n}\over 1-x_{i+n}}-{x_{i+n}\over x_i-x_{i+n}}\right)\\
&=&-\sum_{1\le k\not=i\le n}\left(b_k{x_k\over x_i-x_k}+b_{k+n}{-x_k\over x_i+x_k}\right)\\
&+&{(a+b_i+b_{i+n})x_i^2
+(a-c-b_{2n+1}+1+b_i-b_{i+n})x_i+b_{2n+1}-c+1\over(1-x_i)(1+x_i)},\\
p^{i+n}_{i,i}&+&2p^{i+n}_{i,i+n}+p^{i+n}_{i+n,i+n}=
b_1\left\{{x_i\over1-x_i}-{x_i\over-x_i-x_i}\right\}+2b_i{x_i\over-x_i-x_i}\\
&+&{1\over1+x_i}\left\{(a+b_{i+n})(-x_i)-c+1{-1\over x_i}\right\}
-\sum_{1\le k\not=i+n\le 2n+1}b_k{x_k\over x_{i+n}-x_k}\\
&=&-\sum_{1\le k\not=i\le n}\left(b_k{x_k\over -x_i-x_k}+b_{k+n}{-x_k\over -x_i+x_k}\right)\\
&+&
{(a+b_i+b_{i+n})x_i^2+(-a+c+b_{2n+1}-1+b_i-b_{i+n})x_i+b_{2n+1}-c+1\over(1-x_i)(1+x_i)}.
\end{eqnarray*}
In order the system for $v$ to be of rank $n+1$, $D_iD_jv$ must be linearly 
related with $D_kv$ and $v$. We have
$$p^{2n+1}_{i,i}+p^{2n+1}_{i+n,i+n}=0\quad\mbox{and}\quad p^k_{i,j}+p^k_{i+n,j+n}=p^{k+n}_{i,j}+p^{k+n}_{i+n,j+n}\quad (1\le i,j\le n),$$
if and only if $b_i=b_{i+n}\ (1\le i\le n)$. Assuming these, we
have
$$
p^i_{i,i}+2p^i_{i,i+n}+p^i_{i+n,i+n}=p^{i+n}_{i,i}+2p^{i+n}_{i,i+n}+p^{i+n}_{i+n,i+n}\quad
(1\le i\le n),$$
if and only if $a-c-b_{2n+1}+1=0$.
\begin{prp}The pull-back, under the embedding $\iota:X(2,n+3)\to X(2,2n+4)$ of 
a(ny) non-zero solution of the hypergeometic system $E_D^{2n+1}(a,b_1,\dots,b_{2n+1},c)$
satisfies a system of rank $n+1$ if and only if 
$$b_i=b_{i+n},\ (1\le i\le n)\quad\mbox{and}\quad a-c-b_{2n+1}+1=0.$$
This system of rank $n+1$ on $X(2,n+3)$, which will be called
 ${\mathcal E}(a,b_1,\dots,b_n,c)$, is given by
$$D_iD_jv=\sum_{k=1}^nq_{ij}^kD_kv + q^0_{ij}v\quad (1\le i,j\le n),$$
with parameters $(a, b_1,\dots,b_n,c)$, where
\begin{eqnarray*}
q^k_{ij}&=&q^k_{ji},\quad q^0_{ij}=q^0_{ji},\quad 1\le i,j,k\le n;\quad
q^0_{ij}=q^k_{ij} =0,\quad i\not=j\not=k\not=i\\
q^i_{ij}&=&p^i_{i,j}+p^i_{i,j+n}
=2b_j{x_j^2\over x_i^2-x_j^2},\quad i\not=j\\
q^0_{ii}&=&p^0_{i,i}+p^0_{i+n,i+n}
=a(2b_i){x_i^2\over 1-x_i^2},\\
q^k_{ii}&=&p^k_{i,i}+p^k_{i+n,i+n}
=2b_i\left({x_i^2\over1-x_i^2}-{x_i^2\over x_k^2-x_i^2}\right),\quad i\not=k\\
q^i_{ii}&=&p^i_{i,i}+2p^i_{i,i+n}+p^i_{i+n,i+n}\\
&=&-\sum_{k\not=i}2b_k {x_k^2\over x_i^2-x_k^2} + 
{(a+2b_i)x_i^2 +a- 2c+2\over 1-x_i^2}.
\end{eqnarray*}
The system has regular singularities along 
$$\prod_{i=1}^nx_i(1-x_i)(1+x_i)\prod_{1\le i<j\le n}(x_i-x_j)(x_i+x_j)=0$$
and at infinity.
\end{prp}
Rewriting this system in the form $\partial^2v/\partial x_i\partial x_j=\cdots$, we get
\begin{cor}This system is non-singular along $\{x_j=0\}$ if and only if
$$(q^i_{ii}|_{x_i=0}=)\quad 2\sum_{1\le k(\not=j)\le n}b_k +a-2c+2=0.$$
\end{cor}
We need not check that this system is of rank $n+1$, because of the
 following fact. Note that if we introduce the variables $y_j=x_j^2$,
 then we have
$$D_j^x=2D_j^y,\quad\mbox{where}\quad D_j^x=x_j{\partial/\partial
 x_j},\ D_j^y=y_j{\partial/\partial y_j}.$$
Comparing the coefficients of ${\mathcal E}(a,b,c)$ and those of $E_D^n(a,b,c)$,
 we have
\begin{cor}The Appell-Lauricella system $E_D^n(a/2,b_1,\dots,b_n,c-a/2)$ 
 in $y$-variables is transformed into ${\mathcal E}(a,b,c)$ in $x$-variables by
 the change $y_j=x_j^2$.\end{cor}
\begin{rmk}The following integral representation of a solution of $E_D^{2n+1}(a,b,c)$
$$\int t^{a-1}(1-t)^{c-a-1}(1-x_1)^{-b_1}\cdots(1-x_{2n})^{-b_{2n}}(1-x_{2n+1})^{-b_{2n+1}}dt$$
supports  Proposition 1 and Corollary 2.\end{rmk}
\section{Systems invariant under the involution $\#$}
Conisder the involution 
$$\#:(x_1,\dots,x_n)\longmapsto ({1/ x_1},\dots {1/ x_n})$$
on $X(2,n+3)$ and let $\#$ takes the system $E(a,b,c)=E_D^n(a,b_1,\dots,b_n,c)$
 into a system 
$${}^\#E(a,b,c):D_iD_ju=\sum_{k=1}^n{\pi}_{ij}^kD_ku+{\pi}_{ij}^0u,\quad (1\le
i,j\le n).$$
Since the change of variables $x_i\to 1/x_i$ induces the change of the
derivations $D_i\to -D_i$, we can easily see the coefficients:
\begin{eqnarray*}
{\pi}^k_{ij}&=&{\pi}^k_{ji},\quad {\pi}^0_{ij}={\pi}^0_{ji},\quad 1\le i,j,k\le n;\quad
{\pi}^0_{ij}={\pi}^k_{ij} =0,\quad i\not=j\not=k\not=i,\\ 
{\pi}^i_{ij}&=&b_j{x_i\over x_i-x_j},\quad i\not=j,\qquad
{\pi}^0_{ii}=ab_i{1\over 1-x_i},\\
{\pi}^k_{ii}&=&b_i\left({1\over1-x_i}-{x_k\over x_k-x_i}\right),\quad i\not=k\\
{\pi}^i_{ii}&=&-\sum_{k\not=i}b_k {x_i\over x_i-x_k} + 
{(-c+1)x_i +a+b_i\over 1-x_i}.
\end{eqnarray*}
In particular $E(a,b,c)$ is not invariant under the involution $\#$ for
any choice of parameters not simultaneously zero.
\par\smallskip
Let us transform the systems $E$ and ${}^\#E$ into the normal
forms defined below. We change the unknown $v$ of the system $E(a,b,c)$  as $v=\rho w$,
where $\rho$ is a non-zero function. If we write the new system as
$${}^NE(a,b,c):D_iD_jw=\sum_{k=1}^nP_{ij}^kD_kw+P_{ij}^0w,\quad (1\le i,j\le n)$$
then the coefficients $P_{ij}^k$ are given by
$$P_{ij}^k=p_{ij}^k-{D_i\rho\over \rho}\delta_j^k-{D_j\rho\over \rho}\delta_i^k,$$
where $\delta$ denotes the Kronecker symbol.
Now choose $\rho$ so that the system ${}^NE$ be of {\it normal form}, which means by 
definition,
$$\sum_{k=1}^nP_{kj}^k=0,\quad (1\le j\le n)$$ that is,
$$\sum_{k=1}^np_{kj}^k-(n+1){D_j\rho\over \rho}=0;\quad (1\le j\le n)$$
it is known that if the given system $E$ is integrable then there is a
non-zero function $\rho$ solving the above first order system of
differential equations. Thus we have
$$P_{ij}^k=p_{ij}^k-{\delta_j^k\over n+1}p_i-{\delta_i^k\over n+1}p_j,$$
where 
$$
p_j=\sum_{k=1}^np_{kj}^k=\sum_{k\not=j}p_{kj}^k + p_{jj}^j=\sum_{k\not=j}{b_jx_j+b_kx_k\over x_k-x_j}
+{(a+b_j)x_j-c+1\over 1-x_j}.
$$
The coefficients of the normal form ${}^NE(a,b,c)$ are given by
\begin{eqnarray*}
P_{ii}^k&=&p_{ii}^k={b_ix_i(x_k-1)\over(1-x_i)(x_k-x_i)},\quad i\not=k,\\
P_{ij}^i&=&p_{ij}^i-{p_j\over n+1}={b_jx_j\over x_i-x_j}-{p_j\over n+1},\quad i\not=j\\
P_{ii}^i&=&p_{ii}^i-2{p_j\over n+1}.
\end{eqnarray*}

We next find the normal form 
$${}^{N\#}E(a,b,c):D_iD_jw=\sum_{k=1}^n{\Pi}_{ij}^kD_kw+{\Pi}_{ij}^0w,\quad (1\le
i,j\le n)$$
of ${}^\#E(a,b,c)$. Its coefficients are given by
$${\Pi}_{ij}^k={\pi}_{ij}^k-{\delta_j^k\over n+1}{\pi}_i-{\delta_i^k\over n+1}{\pi}_j,$$
where 
$$
{\pi}_j=\sum_{k=1}^n{\pi}_{kj}^k=\sum_{k\not=j}{\pi}_{kj}^k + {\pi}_{jj}^j=
\sum_{k\not=j}{b_jx_k+b_kx_j\over x_k-x_j}
+{(-c+1)x_j+a+b_j\over 1-x_j}.
$$
Thus we have
\begin{eqnarray*}
{\Pi}_{ii}^k&=&{\pi}_{ii}^k={b_ix_i(x_k-1)\over(1-x_i)(x_k-x_i)},\quad i\not=k,\\
{\Pi}_{ij}^i&=&{\pi}_{ij}^i-{{\pi}_j\over n+1}={b_jx_i\over x_i-x_j}-{{\pi}_j\over n+1},\quad i\not=j\\
{\Pi}_{ii}^i&=&{\pi}_{ii}^i-2{{\pi}_j\over n+1}.
\end{eqnarray*}

Compare the coefficients and recall an elementary fact that a rational function 
$$\sum_j{\alpha_jx+\beta_j\over x-t_j}$$
in $x$ vanishes identically if and only if
$\alpha_jt_j\sum_j+\beta_j=0$ and $\sum_j\alpha_j=0.$
Then we get
\begin{prp}
The system ${}^NE(a,b,c)$ coincides with the system ${}^{N\#}E(a,b,c)$ if and
only if
$$b_1=\cdots=b_n (=b),\quad -nb+a+c-1=0.$$
\end{prp}
Let us summarize the above computaion.
We got a system ${\mathcal E}(a,b,c)={\mathcal E}(a,b_1,\dots,b_n,c)$ of rank $n+1$ on $X(2,n+3)$,
which was the pull-back under $\iota:X(2,n+3)\to X(2,2n+4)$ of the hypergeometric system $E_D^{2n+1}(a,b_1,\dots,b_{2n+1},c)$, with
$$b_{n+1}=b_1,\dots,b_{2n}=b_n,\ b_{2n+1}=a-c+1.$$
The normal forms of ${\mathcal E}(a,b,c)$ and its pull-back coincide
if and only if
$$b_1=\cdots=b_n\ (=b),\quad -(2n+1)b+a+c-1.$$
These conditions lead to 
$$a=(n+1)b,\quad c=nb+1.$$
If moreover ${\mathcal E}$ is regular along the divisors $\{x_j=0\}$, 
then we have $b=1/(n-1)$.
\section{The double covering $f:X(2,7)\to X(3,6)$}
Let $X(3,6)$ be the configuration space of colored six points in general 
position in the projective plane
$$X(3,6)=GL_3\backslash\{z\in M(3,6)\mid \mbox{no 3-minor vanishes}\}
/(\C^\times)^6.$$
We define a rational map $f$ from $X(2,7)$ to $X(3,6)$. 
We start from a system of seven points on the line representing a 
point of $X(2,7)$.
We regard the line, carrying the seven points, as a non-singular conic 
in the plane. 
The five points represented by
the last five points, and the intersection point of the tangent
lines (to the conic) at the first and the second points define a system of 
six points on the
plane, representing a point of $X(3,6)$. Let us express this map $f$
in terms of coordinates. We normalize the system of
seven points to be
$$x=(0,\infty,1,x_1,\dots,x_4)\in X(2,7).$$
If the conic is given by $t_1^2-t_0t_2=0$ in the plane coordinatized by
$t_0:t_1:t_2$, the seven points are represented by the seven columns
$$\begin{array}{ccccccc}
1&0&1&1    &1    &1    &1\\ 
0&0&1&x_1  &x_2  &x_3  &x_4\\ 
0&1&1&x_1^2&x_2^2&x_3^2&x_4^2
\end{array}
.$$
Since the tangents of the conic at the first two points are $t_2=0$ and
$t_0=0$, the intersection is given by $0:1:0$; so the point $f(x)$ is
represented by
$$\left(\begin{array}{cccccc}
0&1&1    &1    &1    &1\\ 
1&1&x_1  &x_2  &x_3  &x_4\\ 
0&1&x_1^2&x_2^2&x_3^2&x_4^2
\end{array}\right)\in X(3,6)
.$$
Normalizing the above $3\times6$-matrix into 
\[ 
\left( 
\begin{array}{cccccc}
1&0&0&1& 1& 1\\ 0&1&0&1&z_1&z_2\\ 0&0&1&1&z_3&z_4
\end{array}
\right),
\]
we can readily show
\begin{prp} The map
$$f:X(2,7)\ni x=(x_1,\dots,x_4)\longmapsto z=(z_1,\dots,z_4)\in X(3,6)$$\end{prp}
defined in this section is given by
\begin{eqnarray*}
z_1&=&{(x_3+x_1)(x_2-1)\over(x_3-1)(x_2+x_1)},\quad z_2={(x_4+x_1)(x_2-1)\over(x_4-1)(x_2+x_1)},\\
z_3&=&{(x_3+1)(x_2-x_1)\over(x_3-x_1)(x_2+1)},\quad z_4={(x_4+1)(x_2-x_1)\over(x_4-x_1)(x_2+1)}.
\end{eqnarray*}
\begin{rmk}\noindent
1.  The Jacobian of $f$ is given by
$${(x_1-1)(-x_2+x_1)(x_2-x_4)(x_3-x_4)(x_3-x_2)(x_1+1)^3(x_2-1)\over
(x_3-x_1)^2(x_2+1)^3(x_4-1)^2(x_2+x_1)^3(x_3-1)^2(-x_4+x_1)^2}. $$
\par\smallskip\noindent
2. Put
\begin{eqnarray*} 
D_1&=&z_1z_4-z_2z_3,\quad D_2=z_1z_4-z_2z_3-z_4+z_2+z_3-z_1,\\
Q&=&-z_2z_3z_1-z_2z_3z_4+z_2z_3+z_1z_4z_2+z_1z_4z_3-z_1z_4.
\end{eqnarray*}
Then the singular locus of the system found in \cite{SaYo2} is defined by
$$D=\prod_{j=1}^4z_j(1-z_j)\cdot(z_1-z_2)(z_1-z_3)(z_2-z_4)(z_3-z_4)D_1D_2Q,$$
and $f^*(D)$ is given by
\begin{eqnarray*}
&-&(x_3+x_1)(x_2-1)^5(x_4+x_1)(x_3+1)(-x_2+x_1)^5(x_4+1)(x_1+1)^{13}\\
&\times&(x_3-x_2)^5(x_2-x_4)^5(x_3-x_4)^5(x_1-1)^5(x_3+x_2)(x_2+x_4)(x_3+x_4)\\
&\times&(x_3-1)^{-7}(x_2+x_1)^{-{11}}(x_4-1)^{-7}(x_3-x_1)^{-7}(x_2+1)^{-{11}}(-x_4+x_1)^{-7}.\end{eqnarray*}
\end{rmk}
\begin{prp}
The map $f$ is invariant under the involution $\#:x_i\to1/x_i.$
Moreover $f$ is a two-to-one map.\end{prp} 
Sketch of the proof.
For given $(z_1,\dots,z_4)$, we solve $x_1,\dots,x_4$. We can see that
$x_2$ must be a solution of the quadratic equation
$$C_2Z^2+C_1Z+C_0=0,$$
where
\begin{eqnarray*}
C_0&=&C_2=-(-z_2z_4+z_3z_2z_4-z_3+z_1z_3+z_4-z_1z_3z_4)\\
&\times&(z_3+z_2-z_4-z_1-z_1z_3z_4+z_1z_3z_2-z_2z_4z_1+z_3z_2z_4-2z_2z_3+2z_1z_4),\\
C_1&=&-6z_1z_3z_4+4z_2^2z_4z_1-4z_2^2z_1z_3+4z_3^2z_2^2z_1+4z_4^2z_1^2z_2-4z_1^2z_2z_4+4z_1^2z_3z_2\\
&-&2z_1z_3z_2-2z_2z_4z_1-6z_3z_2z_4-6z_2z_4^2z_1z_3+16z_2z_4z_1z_3-2z_1^2z_3z_4^2z_2\\
&-&2z_3^2z_1z_2-2z_1^2z_3^2z_2-2z_4^2z_2z_1-2z_1^2z_3z_4+2z_1^2z_3^2z_4-2z_3^2z_2^2z_4^2-2z_3z_2^2z_4\\
&-&2z_2^2z_4^2z_1+2z_2^2z_4^2z_3+4z_3^2z_2^2z_4-2z_1^2z_3^2z_4^2+4z_1^2z_3z_4^2-2z_3^2-2z_4^2\\
&-&4z_3^2z_2^2+2z_1z_3-2z_2z_3+4z_4z_3+4z_4z_1^2-4z_4^2z_1^2+2z_2z_4+4z_1z_4^2-2z_1z_4\\
&+&2z_1z_3^2-2z_1^2z_3+4z_3^2z_2-2z_2^2z_4+2z_2z_4^2+2z_1^2z_3^2z_4z_2+4z_3z_2^2+2z_3z_2^2z_4^2z_1\\
&-&2z_3^2z_2^2z_4z_1-2z_2^2z_4z_1z_3+4z_1z_3^2z_4^2z_2-2z_1^2z_3z_2z_4-6z_1z_3^2z_2z_4.\\
\end{eqnarray*}
The discriminant of the quadratic equation with respect to $Z$ is given by
$$16(z_1-1)(z_2-1)(z_3-1)(z_4-1)(z_1-z_2)(z_3-z_4)D_2Q.$$
Other ones $x_1$ and $x_3$ can be expressed as rational 
functions in $x_2$ and $z$:
\begin{eqnarray*}
x_1&=&-x_2(-z_1z_3+z_1z_3z_2+z_2z_4-z_2z_4z_1-z_2+z_1)/(-z_3x_2+z_3+z_2+z_4x_2\\
&-&z_4-z_1-z_1z_3z_4-z_1z_3z_4x_2+z_1z_3z_2-z_2z_4z_1+z_3z_2z_4x_2+z_3z_2z_4\\
&+&x_2z_1z_3-2z_2z_3-x_2z_2z_4+2z_1z_4),\\
x_3&=&-x_2(1-2z_1+z_1z_3-x_2+x_2z_1z_3)(z_3+z_2-z_4-z_1-z_1z_3z_4+z_1z_3z_2\\
&-&z_2z_4z_1+z_3z_2z_4-2z_2z_3+2z_1z_4)/(2x_2z_3^2z_1z_4+z_3+z_2-z_4-z_1\\
&-&6z_1z_3z_4+2z_1z_3z_2-z_2z_4z_1+z_3z_2z_4+2z_2z_4z_1z_3-4z_3^2z_1z_2\\
&+&z_1^2z_3^2z_2+2z_1^2z_3z_4-z_1^2z_3^2z_4-2z_3^2+2z_1z_3-4z_2z_3+2z_4z_3\\
&+&2z_1z_4+z_1z_3^2-z_1^2z_3+4z_3^2z_2-z_1^2z_3z_2z_4+z_1z_3^2z_2z_4\\
&-&2x_2z_3^2z_2z_4-2x_2z_1^2-2x_2z_1^2z_3z_2+2x_2z_3^2+3x_2z_1^2z_3-x_2z_2\\
&+&x_2z_1-3x_2z_2z_4z_1+x_2z_1z_3^2z_2z_4-2x_2z_4z_3-2z_3^2z_2z\\
&+&2z_1z_3^2z_4-z_3x_2+z_4x_2+2x_2z_1z_2+2x_2z_1^2z_2z_4+3z_3z_2z_4x_2\\
&-&x_2z_1^2z_3z_2z_4-x_2z_1^2z_3^2z_4-3x_2z_1z_3^2+x_2z_1^2z_3^2z_2),\\
%
\end{eqnarray*}
and $x_4$ is obtained from the expression of $x_3$ by the exchanges 
$z_1\leftrightarrow z_2$ and $z_3\leftrightarrow z_4$. Therefore the map $f$ is two-to-one. 

\section{A system on $X(3,6)$ induced by ${\mathcal E}$ invariant under $\#$}
We have now the system 
$${\mathcal E}(a,b,c),\quad a=5b,\ c=4b+1,\ b=b_1=\cdots=b_4,$$
with a parameter $b$ on $X(2,7)$ invariant under the involution
$\#$. The push-down $f_*{\mathcal E}(a,b,c)$ of this system is a system defined
on $X(3,6)$. 
Here we confess honestly that we still do not find a way to express
the system $f_*{\mathcal E}(a,b,c)$ in $z$-variables with parameter $b$ in a 
reasonably compact form. Nevertheless when the system ${\mathcal E}(a,b,c)$ 
is non-sigular
along the divisors $\{x_j=0\}$, that is, when $b=1/3$, we can explicitly find 
the coefficients of the system $f_*{\mathcal E}(a,b,c)$ and we obtain
\begin{thm}The system ${\mathcal E}(5/3,1/3,7/3)$, as a system on on $X(3,6)$,
coincides with the system we
found in \cite{SaYo2}, which has singularities along the divsor $\{D=0\}$.
\end{thm}
Make the coordinate change $x\to z$ to transform the system 
${\mathcal E}(a,b,c)$ into the form 
$\partial^2v/\partial z_i\partial z_j=\cdots$ 
and write the coefficients in terms of $x$. 
Though it is possible to express these coefficients in $z$, 
it is much easier to rewrite the {\it coefficients} of the system in \cite{SaYo2} in termns of $x$. 
We then compare the coefficients of the two systems to find that they coincide.



\begin{thebibliography}{MSTY}
\bibitem{HSY}{\sc M. Hara, T. Sasaki and M. Yoshida,} Tensor products of linear differential equations, Funkcial. Ekvac. 32(1989), 455--477.
\bibitem{MT}{\sc K. Matsumoto and T. Terasoma,} Theta constants associated to cubic threefolds, in preparation, 2000.
\bibitem{SaYo2}{\sc T. Sasaki and M. Yoshida,} 
The uniformizing differential equation of the complex hyperbolic
structure on the moduli space of marked cubic surfaces,
Proc. Japan Acad., 75(1999), 129 --133.
\bibitem{Yo2}{\sc M. Yoshida,} {\em Fuchsian Differential equations}, Vieweg
 Verlag, 1987.
\end{thebibliography}
\end{document}